\documentclass[12pt]{amsart}
\usepackage{amssymb}

\newtheorem{thm}{Theorem}[section]
\newtheorem{lem}[thm]{Lemma}
\newtheorem{cor}[thm]{Corollary}
\newtheorem{prop}[thm]{Proposition}
\newtheorem{ex}[thm]{Example}

\newtheorem*{prob*}{Open problem}

\theoremstyle{definition}

\newtheorem{defi}[thm]{Definition}

\theoremstyle{remark}

\newtheorem{rem}[thm]{Remark}
\newtheorem*{rem*}{Remark}

\DeclareMathOperator{\Aff}{Aff}
\DeclareMathOperator{\Aut}{Aut}

\newcommand{\kringel}{\mathbin{\raise1pt\hbox{$\scriptstyle\circ$}}}
\newcommand{\pkt}{\mathbin{\raise0pt\hbox{$\scriptstyle\bullet$}}}

\newcommand{\C}{\mathbb{C}}

\newcommand{\R}{\mathbb{R}}

\newcommand{\ad}{\mathop{\rm ad}}

\newcommand{\End}{\mathop{\rm End}}
\newcommand{\Der}{\mathop{\rm Der}}

\newcommand{\La}{\mathfrak{a}}
\newcommand{\Lb}{\mathfrak{b}}

\newcommand{\Lf}{\mathfrak{f}}
\newcommand{\Lg}{\mathfrak{g}}

\newcommand{\Ln}{\mathfrak{n}}

\newcommand{\Lr}{\mathfrak{r}}

\newcommand{\al}{\alpha}
\newcommand{\be}{\beta}
\newcommand{\ga}{\gamma}
\newcommand{\de}{\delta}

\newcommand{\la}{\lambda}
\newcommand{\om}{\omega}

\newcommand{\Om}{\Omega}

\newcommand{\ra}{\rightarrow}

\renewcommand{\phi}{\varphi}

\begin{document}


\title[LR-algebras]{LR-algebras}

\author[D. Burde]{Dietrich Burde}
\author[K. Dekimpe]{Karel Dekimpe}
\author[S. Deschamps]{Sandra Deschamps}
\address{Fakult\"at f\"ur Mathematik\\
Universit\"at Wien\\
  Nordbergstr. 15\\
  1090 Wien \\
  Austria}
\email{dietrich.burde@univie.ac.at}
\address{Katholieke Universiteit Leuven\\
Campus Kortrijk\\
8500 Kortrijk\\
Belgium}
\email{karel.dekimpe@kuleuven-kortrijk.be}
\email{sandra.deschamps@kuleuven-kortrijk.be}

\date{\today}

\thanks{The first author thanks the KULeuven Campus Kortrijk  for its hospitality and support}
\thanks{The second author expresses his gratitude towards the Erwin Schr\"odinger International
Institute for Mathematical Physics}
\thanks{Research supported by the Research Programme of the Research Foundation-Flanders (FWO): G.0570.06}
\thanks{Research supported by the Research Fund of the Katholieke Universiteit Leuven}

\begin{abstract}
In the study of NIL-affine actions on nilpotent Lie groups we
introduced  so called LR-structures on Lie algebras. The aim of this
paper is to consider the existence question of LR-structures, and to
start a structure theory of LR-algebras. We show that any Lie
algebra admitting an LR-structure is $2$-step solvable. Conversely
we find several classes of $2$-step solvable Lie algebras admitting
an LR-structure, but also classes not admitting such a structure. We
study also ideals in LR-algebras, and classify low-dimensional
LR-algebras over $\R$.
\end{abstract}

\maketitle

\section{Introduction}

LR-algebras and LR-structures on Lie algebras arise in the study of affine
actions on nilpotent Lie groups as follows. Let $N$ be a real, connected and simply connected 
nilpotent Lie group. Denote by  $\Aff (N)=N \rtimes \Aut(N)$ the group of affine
transformations of $N$, acting on $N$ via
\[ \forall m,n\in N,\;\forall \alpha \in \Aut(N): \;\;^{(m,\alpha)}n=m\cdot \alpha(n).\]
Note that for the special case where $N=\R^n$, we obtain the usual group of affine 
transformations $\Aff(\R^n)$ of $n$-dimensional space.
When $N$ is not abelian, we sometimes talk about the NIL-affine group $\Aff(N)$, or NIL-affine 
motions. Recently, there has been a growing interest in those subgroups $G\subseteq \Aff(N)$ 
which act either properly discontinuously (in case $G$ is discrete) or simply transitively 
(in case $G$ is a Lie group)  on $N$ (see for example \cite{BAU}, \cite{DEK}). It is known that
all groups which appear as such a simply transitive NIL-affine group have to be solvable. 
Conversely for any connected and simply connected solvable Lie group $G$, there exists a nilpotent
Lie group $N$ for which one can find an embedding $\rho:G\rightarrow \Aff(N)$ realizing
$G$ as a subgroup of $\Aff(N)$ acting simply transitively on $N$ (see \cite{DEK}). \\
Nevertheless, it is still an open problem to determine for a given 
$G$ all connected, simply connected nilpotent Lie groups $N$, on which $G$ acts simply 
transitively via NIL-affine motions. Even for the case $G=\R^n$ the problem is non-trivial 
and interesting.
For this case we were able to translate this question in \cite{BDD} to the existence problem 
of an LR-structure on the Lie algebra $\Ln$ of $N$. Indeed, we showed the following result 
(for the definition of a complete LR-structure see $\ref{deflr}$). 

\begin{thm}\label{thm-LR}\cite[Theorem 5.1]{BDD}
Let $N$ be a connected and simply connected nilpotent Lie group of dimension $n$. Then there
exists a simply transitive NIL-affine action of $\R^n$ on $N$ if and only if
the Lie algebra $\Ln$ of $N$ admits a complete LR-structure.
\end{thm}

The aim of this paper is to begin a study of LR-algebras and LR-structures on
Lie algebras. Although LR-algebras arose, as we just explained, in the 
context of Lie algebras over the field $\R$, we will work over
an arbitrary field $k$  of characteristic zero.

\begin{defi}\label{deflr}
An algebra $(A,\cdot)$ over $k$ with product $(x,y) \mapsto x\cdot y$
is called an {\it LR-algebra}, if the product satisfies
the identities
\begin{align}
x\cdot (y\cdot z)& = y\cdot (x\cdot z) \label{lr1}\\
(x\cdot y)\cdot z& =(x\cdot z)\cdot y \label{lr2}
\end{align}
for all $x,y,z \in A$.
\end{defi}
Denote by $L(x), R(x)$ the left respectively right multiplication operator in
the algebra $(A,\cdot)$. The letters LR stand for ``left and right'', indicating
that in an LR-algebra the left and right multiplication operators commute:
\begin{align}
[L(x),L(y)] & = 0, \label{lr3}\\
[R(x),R(y)] & = 0  \label{lr4}.
\end{align}

LR-algebras are Lie-admissible algebras:

\begin{lem}
The commutator $[x,y]=x\cdot y-y\cdot x$ in an LR-algebra $(A,\cdot)$
defines a Lie bracket.
\end{lem}

\begin{proof}
We have, using the above identities for all $x,y,z\in A$,
\begin{align*}
[[x,y],z]+[[y,z],x]+[[z,x],y] & = [x,y]\cdot z - z\cdot [x,y] +
[y,z]\cdot x - x \cdot [y,z] \\
 & + [z,x]\cdot y - y \cdot [z,x] \\
 & = 0.
\end{align*}
This shows that the Jacobi identity is indeed satisfied.
\end{proof}

The associated Lie algebra $\Lg$ then is said to admit
an LR-structure:

\begin{defi}\label{affine}
An {\it LR-structure} on a Lie algebra
$\Lg$ over $k$ is an LR-algebra product $\Lg \times \Lg \rightarrow \Lg$
satisfying
\begin{align}\label{lr5}
[x,y]& =x\cdot y -y\cdot x
\end{align}
for all $x,y,z \in \Lg$. The LR-structure, resp.\ the LR-algebra is said to be
{\it complete}, if all left multiplications $L(x)$ are nilpotent.
\end{defi}

\begin{rem}
If $\Lg$ is abelian, then an LR-structure on $\Lg$ is commutative and associative.
Indeed, then we have $R(x)=L(x)$. Conversely, commutative, associative algebras
form a subclass of LR-algebras with abelian associated Lie algebra.
\end{rem}

To conclude this introduction, let us present some easy examples of LR-algebras.
Denote by $\Lr_2(k)$ the $2$-dimensional non-abelian Lie algebra over $k$ with basis 
$(e_1,e_2)$, and $[e_1,e_2]=e_1$.

\begin{ex}
The classification of non-isomorphic LR-algebras $A$ with associated Lie algebra $\Lr_2(k)$
is given as follows:
\vspace*{0.5cm}
\begin{center}
\begin{tabular}{c|c}
$A$ & Products \\
\hline
$A_1$ & $e_1\cdot e_1=e_1,\; e_2\cdot e_1=- e_1.$\\
\hline
$A_2$ & $e_1\cdot e_2=e_1.$\\
\hline
$A_3$ & $e_2\cdot e_1=-e_1.$\\
\hline
\end{tabular}
\end{center}
\end{ex}
\vspace*{0.5cm}
The proof consists of an easy computation. The left multiplications defining
an LR-algebra with associated Lie algebra $\Lr_2(k)$ are of the following form:
\[
L(e_1)=\begin{pmatrix} \al & \be \\ 0 & 0 \end{pmatrix},\;
L(e_2)=\begin{pmatrix} \be-1 & \ga \\ 0 & 0 \end{pmatrix},
\]
where $\al\ga=\be(\be-1)$. All these algebras are isomorphic to one
of the algebras $A_1,A_2,A_3$. Note that the algebra $A_2$ is
complete, whereas the algebras $A_1$ and $A_3$ are incomplete.

\section{Structural properties of LR-algebras}

We just saw examples of LR-structures on a $2$-step solvable Lie algebra, i.e., on
$\Lr_2(k)$. It turns out that {\it all} Lie algebras admitting an LR-structure are
two-step solvable.

\begin{prop}
Any Lie algebra over $k$ admitting an LR-structure is two-step solvable.
\end{prop}

\begin{proof}
For any $x,y,u,v\in A$ we have the following symmetry relation:
\begin{eqnarray*}
(x\cdot y)(u\cdot v) & = & (x\cdot (u\cdot v))\cdot y \\
         & = & (u\cdot(x\cdot v))\cdot y \\
         & = & (u\cdot y)\cdot (x\cdot v) \\
         & = & x\cdot ((u\cdot y)\cdot v) \\
         & = & x\cdot ((u\cdot v)\cdot y)\\
         & = & (u\cdot v)\cdot (x\cdot y).
\end{eqnarray*}
Using this we obtain
\begin{eqnarray*}
[[x,y],[u,v]] & = & [x\cdot y-y\cdot x,u\cdot v-v\cdot u] \\
& = & (x\cdot y-y\cdot x)\cdot (u\cdot v-v\cdot u) - (u\cdot v-v\cdot u)
\cdot (x\cdot y-y\cdot x) \\
& = & (x\cdot y)\cdot (u\cdot v) - (x\cdot y)\cdot (v\cdot u) -(y\cdot x)
\cdot (u\cdot v) + (y\cdot x)\cdot (v\cdot u) \\
& -  & (u\cdot v)\cdot (x\cdot y) + (v\cdot u)\cdot (x\cdot y)+(u\cdot v)
\cdot (y\cdot x) -(v\cdot u)\cdot (y\cdot x) \\
& = & 0.
\end{eqnarray*}
This shows that the associated Lie algebra is two-step solvable.
\end{proof}

When we translate this result, using Theorem~\ref{thm-LR}, in terms of NIL-affine
actions, we find the following:

\begin{thm}
Let $N$ be a connected and simply connected nilpotent Lie group for which
$\Aff(N)$ contains an abelian Lie subgroup acting simply transitively on
$N$, then $N$ is two-step solvable.
\end{thm}

\begin{rem}
This result also explains Proposition 4.2 of \cite{BDD} in a much more conceptual way.
\end{rem}

We now present some identities, which are  useful when
constructing LR-structures on a given Lie algebra.
The first pair of identities remind one of the Jacobi identity for Lie
algebras:

\begin{lem}\label{lem-jacob}
Let $(A,\cdot)$ be an LR-algebra. For all $x,y,z \in A$ we have:
\begin{align}
[x,y]\cdot z + [y,z]\cdot x + [z,x]\cdot y & = 0, \label{lr6} \\
x \cdot [y,z] + y\cdot [z,x] + z\cdot [x,y] & = 0.\label{lr7}
\end{align}
\end{lem}

\begin{proof}
The first identity holds because we have
\begin{align*}
[x,y]\cdot z + [y,z]\cdot x + [z,x]\cdot y & =  (x\cdot y-y\cdot x)\cdot z
+(y\cdot z-z\cdot y)\cdot x +(z\cdot x-x\cdot z)\cdot y\\
 & = ((x\cdot y)\cdot z -  (x\cdot z)\cdot y) +
     ((y\cdot z)\cdot x -  (y\cdot x)\cdot z) \\
 & \; + ((z\cdot x)\cdot y - (z\cdot y)\cdot x)\\
 & = 0.
\end{align*}
The second identity follows similarly.
\end{proof}

We also have the following operator identities:

\begin{lem}\label{lem-linear}
In an LR-algebra we have the following operator identities:
\begin{align}
\ad([x,y])-[\ad(x),L(y)]-[L(x),\ad(y)] & = 0.\label{lr8} \\
\ad([x,y])+[\ad(x),R(y)]+[R(x),\ad(y)] & = 0. \label{lr9}
\end{align}
\end{lem}
\begin{proof}
Using $\ad (x)=L(x)-R(x)$
and (\ref{lr3}) and (\ref{lr4})
we obtain
\begin{align*}
\ad ([x,y]) & = [\ad (x),\ad (y)] \\
            & = [L(x)-R(x),L(y)-R(y)]\\
            & = [L(x),L(y)]-[R(x),L(y)]-[L(x),R(y)]+[R(x),R(y)]\\
            & = [-R(x),L(y)]+[L(x),-R(y)] \\
            & = [\ad(x),L(y)]+[L(x),\ad(y)]
\end{align*}
This shows the first identity. The second identity follows similarly.
\end{proof}

We now study ideals of LR--algebras. 
\begin{lem}
Let $(A,\cdot)$ be an LR-algebra and $I,J$ be two-sided ideals of $A$.
Then $I \cdot J$ is also a two-sided ideal of $A$.
\end{lem}

\begin{proof}
It is enough to show that for all $a\in A$, $x\in I$ and $y\in J$,
both $a\cdot(x\cdot y)$ and $(x\cdot y) \cdot a $ belong to $I\cdot J$.
But this is easy to see:
\begin{align*}
a\cdot(x\cdot y) & = x \cdot (a \cdot y) \in I \cdot J,\\
(x\cdot y) \cdot a & = (x \cdot a ) \cdot y \in I \cdot J.
\end{align*}
\end{proof}

Before continuing the study of ideals let us note the following:

\begin{lem}\label{derivations}
Let $(A,\cdot)$ be an LR-algebra with associated Lie algebra $\Lg$, and $a\in A$.
Then all operators $L(a)$ and $R(a)$ are Lie derivations of $\Lg$, i.e.,
for any $x,y\in A$, the following identities hold:
\begin{align*}
a\cdot [x,y] & = [a\cdot x, y] + [ x , a\cdot y], \\
[x,y]\cdot a & = [x\cdot a , y] + [x , y \cdot a].
\end{align*}
\end{lem}
\begin{proof}
We have
\begin{align*}
a\cdot [x,y] & = a \cdot (x \cdot y)  - a \cdot (y \cdot x)\\
& = x \cdot ( a \cdot y) - y \cdot (a \cdot x) - (a \cdot y ) \cdot x + (a \cdot x) \cdot y \\
& = [a \cdot x, y] + [x, a \cdot y ].
\end{align*}
The second identity follows similarly.
\end{proof}

The above lemma implies the following result:

\begin{cor}
Let $(A, \cdot)$ be an LR-algebra and assume that $I,J$ are
two-sided ideals of $A$. Then $[I,J]$ is also a two-sided ideal of $A$.
\end{cor}

In particular, $[A,A]$ is a two-sided ideal in $A$.
Let $\ga_1(A)=A$ and $\gamma_{i+1}(A)=[A,\gamma_{i}(A)]$ for all $i\ge 1$.

\begin{cor}\label{ideals}
Let $A$ be an LR-algebra. Then all $\gamma_i(A)$ 
are two-sided ideals of $A$.
\end{cor}

\begin{lem}
Let $A$ be an LR-algebra. Then we have
\[
\gamma_{i+1}(A) \cdot \gamma_{j+1}(A) \subseteq \gamma_{i+j+1}(A)
\]
for all $i,j\ge 0$.
\end{lem}
\begin{proof}
We will use induction on $i\ge 0$. The case $i=0$ follows
from the fact that $\gamma_{j+1}(A)$ is an ideal of $A$. Now assume $i\geq 1$ and
$\gamma_k(A)\cdot \gamma_{j+1}(A)\subseteq \gamma_{k+j}(A)$
for all $k=1,\ldots ,i $.\\
Let $x\in \gamma_1(A)$, $y\in \gamma_i(A)$ and $z\in \gamma_{j+1}(A)$. We have to show
that $[x,y]\cdot z \in \gamma_{i+j+1}(A)$. Using Lemma~\ref{derivations} and
the induction hypothesis, we see that
\[ [x,y]\cdot z = [x \cdot z, y ] + [ x, y \cdot z] \in \gamma_{i+j+1}(A).\]
\end{proof}

It is natural to introduce the center of an LR-algebra $(A,\cdot)$ by
\[
Z(A)=\{ x\in A\mid x\cdot y =y\cdot x \mbox{ for all } y \in A\}.
\]
Clearly $Z(A)$ coincides with $Z(\Lg)$, the center of the associated Lie algebra $\Lg$.

\begin{lem}\label{centercomm}
Let $(A, \cdot)$ be an LR-algebra. Then $Z(A) \cdot [A,A] =[A,A]\cdot Z(A)=0$.
\end{lem}

\begin{proof}
Let $a,b \in A$ and $z \in Z(A)$. By $\eqref{lr7}$ we have
\[
z \cdot [a,b] + a \cdot [b,z] + b \cdot [z,a] = 0.
\]
Since $z\in Z(\Lg)$, where $\Lg$ is the associated Lie algebra
of $A$, we obtain $z \cdot [a,b]=0$. Analogously one shows that $[a,b]\cdot z =0$.
\end{proof}

\begin{lem}\label{centerideal}
Let $A$ be an LR-algebra. Then $Z(A)$ is a two-sided ideal
of $A$.
\end{lem}

\begin{proof}
Let $z \in Z(A)$. We have to show that $[a\cdot z, b ] = [ z\cdot a, b ]=0$
for all $a, b \in A$. Using Lemma~\ref{derivations} we see that
\begin{align*}
a \cdot [z,b] & = [a \cdot z , b] + [z, a \cdot b] \\
[z,b] \cdot a & = [z\cdot a, b ] + [z, b \cdot a].
\end{align*}
Since $z\in Z(A)$ the claim follows.
\end{proof}

Let $Z_1(A)=Z(A)$ and define $Z_{i+1}(A)$ by the identity
$Z_{i+1}(A)/Z_{i}(A)= Z(A/Z_i(A))$. Note that the $Z_i(A)$ are the terms of the
upper central series of the associated Lie algebra $\Lg$.
As an immediate consequence of the previous lemma, we obtain

\begin{cor}
Let $A$ be an LR-algebra. Then all $Z_i(A)$ are two-sided ideals of $A$.
\end{cor}

\section{Classification of LR-structures}

A classification of LR-structures in general is as hopeless as a classification of
Lie algebras. However, one can study such structures in low dimensions.
We will give here a classification of complete LR-structures on real nilpotent
Lie algebras of dimension $n\le 4$. The restriction to complete structures reduces
the computations a lot, in particular for abelian Lie algebras.
Nevertheless, we have classified also incomplete LR-structures in some cases.

If the Lie algebra is abelian then LR-structures, and also LSA-structures, are just 
given by commutative and associative algebras. Here a classification in terms of
polynomial rings and their quotients is well known for $n\le 6$, see \cite{PON} 
and the references cited therein. We would like, however, to have explicit lists in terms
of algebra products. This seems only available in dimension $n\le 3$ 
over $\R$ and $\C$, see \cite{GOR}. 
For $n=4$, there is an explicit list of {\it nilpotent} commutative, associative
algebras (see the references in \cite{PON}), but not for all algebras. 
Such nilpotent, commutative, associative algebras correspond exactly to complete 
left-symmetric algebras with abelian associated Lie algebra.  As \cite{KIM} gives a list 
of all complete LSAs with a nilpotent associated Lie algebra in dimension $4$, 
one can easily extract those with an abelian associated Lie algebra from that list, and 
so one obtains the complete abelian LR-structures  in dimension 4.

It remains to classify all complete LR-structures on a non-abelian nilpotent
Lie algebra of dimension $n\le 4$ over $\R$, which is one of the following:
\vspace*{0.5cm}
\begin{center}
\begin{tabular}{c|c}
$\Lg$ & Lie brackets \\
\hline
$\Ln_3(\R)$ & $[e_1,e_2]=e_3$ \\
\hline
$\Ln_3(\R) \oplus  \R$ & $[e_1,e_2]=e_3$ \\
\hline
$\Ln_4(\R)$ & $ [e_1,e_2]=e_3, [e_1,e_3]=e_4$ \\
\end{tabular}
\end{center}

\begin{prop}
The classification of LR-algebra structures on the Heisenberg Lie algebra
$\Ln_3(\R)$ is given as follows:
\vspace*{0.5cm}
\begin{center}
\begin{tabular}{c|c}
$A$ & Products \\
\hline
$A_1(\al),\, \al\in \R$ & $e_1\cdot e_1=e_3,\; e_1\cdot e_2=e_3,\; e_2\cdot e_2=\al e_3.$\\
\hline
$A_2(\be),\, \be \in \R$ & $e_1\cdot e_2=\be e_3,\; e_2\cdot e_1=(\be-1)e_3,
\; e_2\cdot e_2=e_1.$ \\
\hline
$A_3$ & $e_1\cdot e_2=\frac{1}{2}e_3,\; e_2\cdot e_1=-\frac{1}{2}e_3.$\\
\hline
$A_4$ & $e_2\cdot e_1=-e_3,\; e_2\cdot e_2=e_2,$\\
      & $e_2\cdot e_3=e_3,\;e_3\cdot e_2=e_3$ \\
\hline
\end{tabular}
\end{center}
\vspace*{0.5cm}
All LR-algebras are complete, except for $A_4$.
\end{prop}

\begin{proof}
Any LR-algebra structure on $\Ln_3(\R)$ is isomorphic to one of the following, written
down with left multiplication operators for the basis $(e_1,e_2,e_3)$:
\[
L(e_1)=\begin{pmatrix} 0 & 0 & 0 \\ \al & \ga & 0 \\ \be & \de & \ga \end{pmatrix},\;
L(e_2)=\begin{pmatrix} 0 & \la & 0 \\ \ga & \mu & 0 \\ \de-1 & \nu & \mu \end{pmatrix},\;
L(e_3)=\begin{pmatrix} 0 & 0 & 0 \\ 0 & 0 & 0 \\ \ga & \mu & 0 \end{pmatrix},
\]
satisfying the following polynomial equations:
\begin{align*}
\al\la & = 0 \\
\ga\la & = 0 \\
\ga^2-\al\mu & = 0 \\
\ga(2\de-1)-\al\nu-\be\mu & = 0\\
\be\la  & = 0
\end{align*}
A straightforward case by case analysis yields the result. We have
$A_1(\al)\simeq A_1(\al')$ if and only $\al'=\al$, and the same result for
$A_2(\be)$.
\end{proof}

Similarly we obtain the following result.

\begin{prop}
The classification of LR-algebra structures on $\Lg=\Ln_4(\R)$
is given as follows:
\vspace*{0.5cm}
\begin{center}
\begin{tabular}{c|c}
$A$ & Products \\
\hline
$A_1(\al)$ & $e_1\cdot e_1=\al (\al-1)e_2,\; e_1\cdot e_2=\al e_3,\; e_1\cdot e_3=\al e_4,$\\
$\al\in \R$  & $e_2\cdot e_1=(\al-1)e_3,\; e_2\cdot e_2=e_4,\;e_3\cdot e_1=(\al-1)e_4.$ \\
\hline
$A_2$ & $e_1\cdot e_1=e_3,\; e_2\cdot e_1=-e_3,$\\
           & $e_2\cdot e_2=e_4,\;e_3\cdot e_1=-e_4.$ \\
\hline
$A_3$ & $e_1\cdot e_1=e_3,\; e_1\cdot e_2=e_3,$\\
           & $e_1\cdot e_3=e_4,\;e_2\cdot e_2=e_4.$ \\
\hline
$A_4(\al,\be,\ga)$ & $e_1\cdot e_1=\al e_2,\; e_1\cdot e_2=\be e_3+\ga e_4,\;
e_1\cdot e_3=\be e_4,$ \\
$\al,\be,\ga\in \{ 0,1 \}$   &  $e_2\cdot e_1=(\be-1)e_3+\ga e_4,\; e_3\cdot e_1=(\be-1)e_4.$ \\
\hline
$A_5(\al)$ & $e_1\cdot e_1=\al e_4,\; e_2\cdot e_1=-e_3,\;e_2\cdot e_2=e_3,$ \\
$\al \in \{ 0,1 \}$   & $e_2\cdot e_3=e_4,\; e_3\cdot e_1=-e_4,\;e_3\cdot e_2=e_4.$ \\
\hline
$A_6$ & $e_2\cdot e_1=- e_3,\; e_2\cdot e_2=e_2,\;e_2\cdot e_3=e_3,\; e_2\cdot e_4=e_4,$ \\
      & $e_3\cdot e_1=-e_4,\;e_3\cdot e_2=e_3,\;e_3\cdot e_3=e_4,\;
e_4\cdot e_2=e_4.$ \\
\hline
\end{tabular}
\end{center}
\vspace*{0.5cm}
The algebra $A_6$ is not complete. All the other ones are complete.
\end{prop}

The family $A_4(\al,\be,\ga)$ consists of $8$ different algebras.

\begin{prop}
The classification of complete LR-algebra structures on $\Lg=\Ln_3(\R)\oplus \R$
is given as follows:
\vspace*{0.5cm}
\begin{center}
\begin{tabular}{c|c}
$A$ & Products \\
\hline
$A_1(\al)$ & $e_1\cdot e_2=\al e_3,\; e_2\cdot e_1=(\al-1) e_3,$\\
$\al\in \R$  & $e_2\cdot e_2= e_1,\; e_4\cdot e_4=e_3.$ \\
\hline
$A_2(\al)$ & $e_1\cdot e_1=\al e_3,\; e_1\cdot e_2=e_4,\; e_2\cdot e_1=-e_3+e_4,$\\
$\al\in  \{ 0,1 \}$   & $e_2\cdot e_2=e_1,\;e_2\cdot e_4=\al e_3,\; e_4\cdot e_2=\al e_3.$ \\
\hline
$A_3(\al,\be)$ & $e_1\cdot e_2=\al e_3,\; e_2\cdot e_1=(\al-1)e_3,\; e_2\cdot e_2=e_1,$\\
$\al\in \R,\; \be\in  \{ 0,1 \}$  & $e_2\cdot e_4=\be e_3,\; e_4\cdot e_2=\be e_3.$ \\
\hline
$A_4(\al)$ & $e_1\cdot e_1=e_4,\; e_1\cdot e_4= e_3,\; e_2\cdot e_1=- e_3,$ \\
$\al\in \{ 0,1 \}$  &  $e_2\cdot e_2=\al e_3,\; e_4\cdot e_1=e_3.$ \\
\hline
$A_5(\al)$ & $e_1\cdot e_4=e_3,\; e_2\cdot e_1=-e_3,$\\
$\al\in \{ 0,1 \}$ & $e_2\cdot e_2=\al e_3,\;e_4\cdot e_1=e_3.$ \\
\hline
\end{tabular}
\end{center}

\begin{center}
\begin{tabular}{c|c}
$A$ & Products \\
\hline
$A_6(\al)$ & $e_1\cdot e_1=\al e_3,\; e_2\cdot e_1=-e_3,$ \\
$\al \in \R$  & $e_2\cdot e_2=e_3,\; e_4\cdot e_4=e_3.$ \\
\hline
$A_7(\al)$ & $e_1\cdot e_1=\al e_3,\; e_2\cdot e_1=-e_3,$ \\
$\al \le \frac{3}{4}$  & $e_2\cdot e_2=-e_3,\;e_4\cdot e_4=e_3.$ \\
\hline
$A_8$ & $e_1\cdot e_2=\frac{1}{2} e_3,\; e_2\cdot e_1=-\frac{1}{2}e_3,$ \\
  & $e_4\cdot e_4=e_3.$ \\
\hline
$A_{9}(\al)$ & $e_1\cdot e_1=e_4,\; e_1\cdot e_2=\al e_3,$ \\
$\al \ge \frac{1}{2}$  & $e_2\cdot e_1=(\al-1)e_3,\;e_2\cdot e_2=e_4.$ \\
\hline
$A_{10}(\al)$ & $e_1\cdot e_1=e_4,\; e_1\cdot e_2=\al e_3,$ \\
$\al \ge \frac{1}{2}$  & $e_2\cdot e_1=(\al-1)e_3,\;e_2\cdot e_2=-e_4.$ \\
\hline
$A_{11}(\al)$ & $e_1\cdot e_1=e_4,\; e_1\cdot e_2=\al e_3,$ \\
$\al \in \R$  & $e_2\cdot e_1=(\al-1)e_3.$ \\
\hline
$A_{12}$ & $e_1\cdot e_1= e_4,\; e_2\cdot e_1=- e_3,$ \\
  & $e_2\cdot e_2=e_3.$ \\
\hline
$A_{13}(\al)$ & $e_1\cdot e_1=e_3,\; e_2\cdot e_1=- e_3,$ \\
$\al \in \R$  & $e_2\cdot e_2=\al e_3.$ \\
\hline
$A_{14}$ & $e_1\cdot e_2= \frac{1}{2}e_3,\; e_2\cdot e_1=-\frac{1}{2}e_3.$ \\
\hline
$A_{15}(\al)$ & $e_1\cdot e_1=e_4,\; e_2\cdot e_1=- e_3,$ \\
$\al \ge 1$  & $e_2\cdot e_2=\al e_3-e_4.$ \\
\hline
\end{tabular}
\end{center}
\end{prop}
\vspace*{0.5cm}
\begin{rem}
The computations for the above result are quite complicated, but do not give
much insight. Therefore we have omitted them here. However, we did the computations
independently to be sure that they are correct.
\end{rem}

\section{LR-structures on nilpotent Lie algebras}

We know that any Lie algebra admitting an LR-structure must be
$2$-step solvable. Conversely we can ask which $2$-step solvable
Lie algebras admit an LR-structure. We start with $2$-step solvable,
filiform nilpotent Lie algebras $\Lf_n$ of dimension $n$.
There exists a so called adapted basis $(e_1,\ldots ,e_n)$ of $\Lf_n$ such that the Lie
brackets are given as follows:
\begin{align*}
[e_1,e_i] & = e_{i+1}, \quad 2\le i\le n-1, \\
[e_2,e_i] & =\sum_{k=i+2}^n c_{i,k}e_k, \quad 3\le i\le n-2, \\
[e_i,e_j] & = 0, \quad 3\le i\le j.
\end{align*}

The Jacobi identity is satisfied if and only if $c_{i+1,k}=c_{i,k-1}$
for all $6\le i+3\le k\le n$. For details, see for example \cite{BRA}.

\begin{lem}
Let $\Lf_n$ be given as above. Then the identities
\begin{align}
\ad(e_2)\ad(e_1)^2 & = \ad (e_1)\ad(e_2)\ad(e_1), \label{fili1}\\
\ad(e_1)\ad(e_2)^2 & = \ad (e_2)\ad(e_1)\ad(e_2), \label{fili2}\\
\ad(e_{i+2}) & = \ad(e_1)^i \ad(e_2)-\ad(e_2)\ad(e_1)^i, \quad i\ge 1. \label{fili3}
\end{align}
hold.
\end{lem}

\begin{proof}
The identity \eqref{fili1} is equivalent to
\begin{align*}
0 & = [\ad(e_1),\ad(e_2)]\ad (e_1) = \ad(e_3)\ad (e_1).
\end{align*}
But this follows from  $[e_3,[e_1,e_k]]=0$ for all $k\ge 1$. Similarly,
\eqref{fili2} is equivalent to $\ad (e_3)\ad (e_2)=0$, which follows again by definition.
The identity \eqref{fili3} is proved by induction on $i\ge 1$.
For $i=1$ we have
\[
\ad(e_3)=[\ad(e_1),\ad(e_2)]=\ad(e_1)\ad(e_2)-\ad(e_2)\ad(e_1).
\]
By induction hypothesis, $\ad (e_{i+1})= \ad(e_1)^{i-1}\ad(e_2)-\ad(e_2)\ad(e_1)^{i-1}$.
Then, using \eqref{fili1} repeatedly, we obtain for $i\ge 2$
\begin{align*}
\ad(e_{i+2}) & = \ad(e_1)\ad(e_{i+1})-\ad(e_{i+1})\ad(e_1) \\
            & =\ad(e_1)^i\ad(e_2)-\ad(e_1)\ad(e_2)\ad(e_1)^{i-1} -\ad(e_1)^{i-1}\ad(e_2)\ad(e_1)
+ \ad(e_2)\ad(e_1)^i \\
            & = \ad(e_1)^i\ad(e_2) -  \ad(e_2)\ad(e_1)^i -  \ad(e_2)\ad(e_1)^i +
\ad(e_2)\ad(e_1)^i \\
            & = \ad(e_1)^i\ad(e_2)-\ad(e_2)\ad(e_1)^i.
\end{align*}
\end{proof}

\begin{prop}
Any $2$-step solvable filiform nilpotent Lie algebra $\Lf_n$ admits a
complete LR-structure.
\end{prop}

\begin{proof}
Define an LR-structure on $\Lf_n$ as follows:
\begin{align*}
L(e_1)& =0, \\
L(e_i) & = \ad(e_1)^{i-2}\ad(e_2), \quad 2\le i\le n.
\end{align*}
In particular, this means
\begin{align*}
e_1\cdot e_j & = 0, \quad e_2\cdot e_j = [e_2,e_j], \quad 1\le j\le n, \\
e_j\cdot e_1 & = [e_j,e_1], \quad e_j\cdot e_2 = 0, \quad 1\le j\le n,
\end{align*}
so that $R(e_1)=-\ad (e_1)$ and $R(e_2)=0$. Furthermore, we have
\begin{align*}
e_i\cdot e_j & = [e_2,e_{i+j-2}], \quad 3\le i,j\le n.
\end{align*}
To see this, note that $e_3\cdot e_j = \ad(e_1)\ad(e_2)(e_j) = [e_2,e_{j+1}]$
for $j\ge 3$. Then the result for $i\ge 3$ follows inductively. \\[0.2cm]
Now let us prove that
\[
e_i\cdot e_j-e_j\cdot e_i=[e_i,e_j], \quad 1\le i\le j\le n.
\]
The cases $i=1$ and $i=2$ are obvious. For $j\ge i\ge 3$ we have
\[
e_i\cdot e_j-e_j\cdot e_i = 0 = [e_i,e_j].
\]
In particular it follows $R(e_i)=L(e_i)-\ad (e_i)$. The formula \eqref{fili3} then implies
\begin{align}
R(e_i)& =\ad (e_2)\ad (e_1)^{i-2}, \quad i\ge 3. \label{fili4}
\end{align}
It remains to show that all operators $L(e_i)$ commute, and all $R(e_i)$ commute, i.e.,
\begin{align*}
L(e_i)L(e_j) & = L(e_j)L(e_i), \quad 1\le i<j \le n, \\
R(e_i)R(e_j) & = R(e_j)R(e_i), \quad 1\le i<j \le n.
\end{align*}
The first identity is obvious for $i=1$. For $2\le i<j \le n$ use \eqref{fili1}
and \eqref{fili2} repeatedly to obtain
\begin{align*}
L(e_i)L(e_j) & = \ad(e_1)^{i-2}\ad(e_2)\ad(e_1)^{j-2}\ad(e_2) \\
             & = \ad(e_1)^{j-2}\ad(e_2)\ad(e_1)^{i-2}\ad(e_2) \\
             & = L(e_j)L(e_i).
\end{align*}
This argument also shows $R(e_i)R(e_j)= R(e_j)R(e_i)$ for $3\le i<j\le n$, because of
\eqref{fili4}. For $i=2$ this is trivially true since $R(e_2)=0$.
For $i=1$ and $j\ge 3$ we have to show that $\ad(e_1)R(e_j)=R(e_j)\ad(e_1)$. This
follows again from \eqref{fili1}.
It is obvious that all $L(e_i)$ are nilpotent, hence the LR-structure is complete.
\end{proof}

It is natural to ask which other nilpotent Lie algebras admit LR-structures.
We first observe the following fact.

\begin{prop}
Every $2$-step nilpotent Lie algebra $\Lg$ admits a complete
LR-structure.
\end{prop}

\begin{proof}
For $x\in \Lg$ define an LR-structure by
\[
L(x)= \frac{1}{2}\ad (x).
\]
Indeed, for all $x,y,z\in \Lg$ we have
\begin{align*}
x\cdot y - y\cdot x & = \frac{1}{2}[x,y] - \frac{1}{2}[y,x] = [x,y], \\
x\cdot (y\cdot z)   & = 0  = y\cdot (x\cdot z),\\
(x\cdot y) \cdot z  & = 0  =  (x\cdot z) \cdot y.
\end{align*}
Finally $L(x)$ is a nilpotent derivation for all $x\in \Lg$, since $\Lg$ is
nilpotent.
\end{proof}

\begin{prop}\label{free}
Every free $3$-step nilpotent Lie algebra $\Lg$ on $n$
generators $x_1,\ldots ,x_n$ admits a complete LR-structure.
\end{prop}

\begin{proof}
A vector space basis of $\Lg$ is given by
\begin{align*}
x_1,x_2, & \cdots ,x_n \\
y_{i,j} & = [x_i,x_j], \quad 1\le i<j\le n,\\
z_{i,j,k} & = [x_i,y_{j,k}].
\end{align*}
An LR-structure on $\Lg$ is defined as follows:
\begin{align*}
x_j\cdot x_i     & =-y_{i,j}, \quad 1\le i<j \le n \\
x_i\cdot y_{j,k} & = z_{i,j,k}, \quad 1\le j<k\le i \\
                 & = z_{k,j,i}, \quad j<i<k \\
y_{j,k}\cdot  x_i & = z_{k,j,i} - z_{i,j,k}, \quad j<i<k \\
                 & =  - z_{i,j,k}, \quad i\le j <k \\
\end{align*}
and all other products equal to zero.
\end{proof}

\begin{ex}
Let $\Lf$ be the free $3$-step nilpotent Lie algebra with $3$ generators.
Then there is a basis $(x_1,\ldots ,x_{14})$ of $\Lf$
with generators $(x_1,x_2,x_3)$ and Lie brackets
\vspace*{0.5cm}
\begin{center}
\parbox{6cm}{$
\begin{array}{rl}
x_4 & = [x_1,x_2] \\
x_5 & = [x_1,x_3] \\
x_6 & = [x_2,x_3] \\
x_7 & = [x_1,[x_1,x_2]] = [x_1,x_4]\\
x_8 & = [x_2,[x_1,x_2]] = [x_2,x_4]\\
x_9 & = [x_3,[x_1,x_2]] = [x_3,x_4]
\end{array}$}
\parbox{6cm}{$
\begin{array}{rl}
x_{10} & = [x_1,[x_1,x_3]] = [x_1,x_5]\\
x_{11} & = [x_2,[x_1,x_3]] = [x_2,x_5]\\
x_{12} & = [x_3,[x_1,x_3]] = [x_3,x_5]\\
x_{11}-x_9 & = [x_1,[x_2,x_3]] = [x_1,x_6]\\
x_{13} & = [x_2,[x_2,x_3]] = [x_2,x_6]\\
x_{14} & = [x_3,[x_2,x_3]] = [x_3,x_6]\\
\end{array}$}
\end{center}
\vspace*{0.5cm}
An LR-structure is given by
\begin{align*}
x_2.x_1 & = -x_4, && x_3.x_6=x_{14}, \\
x_2.x_4 & =  x_8, && x_4.x_1=-x_7, \\
x_2.x_5 & =  x_9, && x_5.x_1=-x_{10}, \\
x_3.x_1 & = -x_5, && x_5.x_2=x_9-x_{11}, \\
x_3.x_2 & = -x_6, && x_6.x_1=x_9-x_{11}, \\
x_3.x_4 & =  x_9, && x_6.x_2=-x_{13}. \\
x_3.x_5 & = x_{12},
\end{align*}
\end{ex}

Proposition $\ref{free}$ implies, in the same way as for Novikov
structures (see \cite{BD}), the following corollary.

\begin{cor}
Any $3$-generated $3$-step nilpotent Lie algebra admits a complete LR-structure.
\end{cor}

One might ask whether or not all $3$-step nilpotent Lie algebras admit an LR-structure.
This turns out to be not the case. To find a counterexample we have to look at
Lie algebras with at least $4$ generators.

\begin{prop}\label{counterexample}
Let $\Lg$ be the following $3$-step nilpotent Lie algebra on $4$ generators
of dimension $13$, with basis $(x_1,\ldots,x_{13})$ and non-trivial Lie brackets
\begin{align*}
[x_1,x_2] & = x_5,   && [x_3,x_4]=-x_5, \\
[x_1,x_4] & = x_6,   && [x_3,x_5]=-x_{11}, \\
[x_1,x_6] & = x_{10}, && [x_3,x_8]=x_{9}, \\
[x_1,x_7] & = x_{11}, && [x_4,x_5]=-x_{12},\\
[x_1,x_8] & = x_{12}, && [x_4,x_6]=x_{9},\\
[x_2,x_3] & = x_7,    &&  [x_4,x_7]=x_9+x_{13}.\\
[x_2,x_4] & = x_8, \\
[x_2,x_5] & = x_{13}, \\
[x_2,x_7]  &= x_{13}, \\
\end{align*}
This $2$-step solvable Lie algebra does not admit an LR-structure.
\end{prop}

\begin{proof}
We will assume that $\Lg$ admits an LR-structure and show that this
leads to a contradiction. We denote by $\ad(x_i)$ the adjoint operators,
by $L(x_i)$ the left multiplications, and by $R(x_i)$ the right multiplication
operators with respect to the basis $(x_1,x_2,\ldots, x_{13})$.
The operators $\ad(x_i)$ are given by the Lie brackets of $\Lg$, while the left
multiplication operators are unknown. We denote the $(j,k)$-th entry of
$L(x_i)$ by
\[
L(x_i)_{j,k} = x^i_{j,k}.
\]
The $j$-th column of $L(x_i)$ gives the coordinates of $L(x_i)(x_j)$.
We have to satisfy the identities \eqref{lr1},
\eqref{lr2} and \eqref{lr5}, where $x$, $y$ and $z$ run over all
basis vectors. This leads to a huge system of quadratic equations in
the variables $x^i_{j,k}$ for $1\leq i,j,k \leq 13$, summing up to a
total of $13^3=2197$ variables. It is quite impossible
to solve these equations without further information. However, we can use
our knowledge on ideals in LR-algebras to conclude that
a lot of unknowns $x^{i}_{j,k}$ already have to be zero. This, 
together with Lemmas \ref{lem-jacob} and \ref{lem-linear},
simplifies the system of equations considerably. In this way we can show
that the equations are contradictory.
This works exactly as in the proof of proposition $3.3$ of our paper
\cite{BDV}, where we proved that the above Lie algebra does not admit a
Novikov structure.
\end{proof}

In \cite{BD} we showed that the $2$-generated, free $4$-step nilpotent
Lie algebra (which is $2$-step solvable) does not admit any Novikov structure.
It turns out however, that this example does admit an LR-structure: \\
Let  $\Lg$ be the free $4$-step nilpotent Lie algebra on $2$ generators $x_1$
and $x_2$. Let $(x_1,\ldots,x_8)$ be a basis  of $\Lg$ with the following Lie brackets:

\begin{center}
\parbox{6cm}{$\begin{array}{rl}
x_3 & = [x_1,x_2] \\
x_4 & = [x_1,[x_1,x_2]] = [x_1, x_3]\\
x_5 & = [x_2,[x_1,x_2]] = [x_2,x_3] \\
 &
\end{array}$}
\parbox{6.5cm}{$\begin{array}{rl}
x_6 & = [x_1,[x_1,[x_1,x_2]]] =  [x_1, x_4] \\
x_7 & = [x_2,[x_1,[x_1,x_2]]] = [x_2,x_4] \\
    & = [x_1,[x_2,[x_1,x_2]]] =  [x_1, x_5] \\
x_8 & = [x_2,[x_2,[x_1,x_2]]] =  [x_2,x_5]
\end{array}$}
\end{center}

\begin{prop}
The  $2$-step solvable Lie algebra from the above example does admit an
LR-structure.
\end{prop}

\begin{proof}
We define the left multilications by $L(x_1)=0$, $L(x_2)=\ad (x_2)$ and,
if $x_i,$ $  i\ge 3$ is a bracket of $x_1$ and $x_2$, then $L(x_i)$ is the
corresponding composition of $\ad (x_1)$ and $\ad (x_2)$:
\begin{align*}
L(x_3) & = L([x_1,x_2])=\ad (x_1)\ad(x_2),\\
L(x_4) & = L([x_1,[x_1,x_2]])=\ad (x_1)^2\ad (x_2)\\
L(x_5) & = L([x_2,[x_1,x_2]])=\ad (x_2)\ad (x_1)\ad (x_2)\\
L(x_6) & = L([x_1,[x_1,[x_1,x_2]]])=\ad(x_1)^3\ad (x_2)\\
L(x_7) & = L([x_2,[x_1,[x_1,x_2]]])=\ad (x_2)\ad (x_1)^2\ad(x_2)\\
       & = L([x_1,[x_2,[x_1,x_2]]])=\ad(x_1)\ad(x_2)\ad(x_1)\ad(x_2)\\
L(x_8) & = L([x_2,[x_2,[x_1,x_2]]])=\ad(x_2)^2\ad(x_1)\ad(x_2).
\end{align*}
In fact this really defines an LR-structure. Note that $L(x_6)=L(x_7)=L(x_8)=0$.
\end{proof}

\section{Construction of LR-structures via extensions}

In the following we will consider Lie algebras $\Lg$ which are an extension
of a Lie algebra $\Lb$ by an abelian Lie algebra $\La$. Hence we have a
short exact sequence of Lie algebras
\begin{equation*}
0 \rightarrow \La \xrightarrow{\iota} \Lg \xrightarrow{\pi}
\Lb \rightarrow 0.
\end{equation*}
Since $\La$ is abelian, there exists a natural $\Lb$-module structure
on $\La$. We denote the action of $\Lb$ on $\La$ by $(x,a)\mapsto \phi (x)a$,
where $\phi \colon \Lb \ra \End (\La)$ is the corresponding Lie algebra
representation. We have
\begin{equation}
\phi([x,y])=\phi(x)\phi(y)-\phi(y)\phi(x)\label{5}
\end{equation}
for all $x,y \in \Lb$.
The extension $\Lg$ is determined by a two-cohomology class. Let
$\Om \in Z^2(\Lb,\La)$ be a $2$-cocycle describing the extension $\Lg$. This implies that
$\Om : \Lb \times \Lb \ra \La$ is a skew-symmetric bilinear map satisfying
\begin{equation}\label{6}
\phi(x)\Om(y,z)-\phi(y)\Om(x,z)+\phi(z)\Om (x,y) =\Om ([x,y],z)-\Om ([x,z],y)+\Om([y,z],x),
\end{equation}
such that the Lie algebra with underlying vector space
$\La \times \Lb$ and Lie bracket given by
 \begin{equation}\label{lie-algebra}
[(a,x),(b,y)]:=(\phi (x)b-\phi(y)a+\Om(x,y),[x,y])
\end{equation}
for $a,b\in \La$ and $x,y\in \Lb$, is isomorphic to $\Lg$.
As a shorthand, we will use $\Lg=(\La,\Lb,\phi,\Om)$ to say that $\Lg$ is the extension
determined by this specific data.

\medskip

Note that the Lie algebras $\Lg$ we are interested in,  are all $2$-step solvable Lie algebras and hence can be
obtained as extensions of two abelian Lie algebras $\La=[\Lg,\Lg]$ and $\Lb=\Lg/[\Lg,\Lg]$. So, although
we will treat extensions in the general case, we
pay special attention to this specific situation where both $\La$ and $\Lb$ are abelian.
In this specific case, the Lie bracket of
$\Lg=\La \times \Lb$ is given by
 \begin{equation*}
[(a,x),(b,y)]:=(\phi (x)b-\phi(y)a+\Om(x,y),0)
\end{equation*}
and the conditions on $\phi$ and $\Om$ are now given as follows: since $\La$ and $\Lb$
are abelian, $\phi$ is just a linear map satisfying
\begin{equation*}
\phi(x)\phi(y)=\phi(y)\phi(x)
\end{equation*}
for all $x,y \in \Lb$.
On the other hand, $\Om : \Lb \times \Lb \ra \La$
is a skew-symmetric bilinear map satisfying
\begin{equation*}
\phi(x)\Om(y,z)-\phi(y)\Om(x,z)+\phi(z)\Om (x,y) =0.
\end{equation*}

Now, let us return to the more general case (i.e., $\Lb$ does not have to be abelian),
and try to construct LR-structures on Lie algebras $\Lg$ which are
given as extension $\Lg=(\La,\Lb,\phi,\Om)$ of a Lie algebra $\Lb$ by an abelian Lie algebra $\La$.
Suppose that we have already an LR-product  $(a,b)\mapsto a\cdot b$ on
$\La$ and an LR-product $(x,y)\mapsto x\cdot y$ on $\Lb$. In other words,
we have
\begin{align*}
x\cdot y-y\cdot x & = [x,y]\\
x\cdot (y\cdot z) & = y\cdot (x\cdot z)\\
(x\cdot y) \cdot z & = (x\cdot z)\cdot y\\
a\cdot b & = b\cdot a \\
a\cdot (b\cdot c) & = b\cdot (a\cdot c)\\
(a\cdot b) \cdot c & = (a\cdot c)\cdot b\\
\end{align*}
for all $x,y,z\in \Lb$ and for all $a,b,c\in \La$. (In fact the product on $\La$ has to be commutative
and associative). We want to lift these LR-products to $\Lg$.
Consider
\begin{align*}
\om & \colon \Lb \times \Lb \ra \La\\
\phi_1,\, \phi_2 & \colon \Lb \ra\End (\La)
\end{align*}
where $\om$ is a bilinear map and
$\phi_1,\, \phi_2$ are linear maps.
We will define a bilinear product $\Lg \times \Lg \ra \Lg$ by
\begin{equation}\label{produkt}
(a,x)\kringel (b,y):=(a\cdot b+\phi_1(y)a+\phi_2 (x)b+\om(x,y),x\cdot y)
\end{equation}

\begin{prop}\label{phi12}
The above product defines an LR-structure on $\Lg$
if and only if the following conditions hold:
\begin{align}
\om(x,y)-\om(y,x) & = \Om(x,y) \label{9}\\
\phi_2(x)-\phi_1(x) & = \phi(x)\label{10}\\
\phi_2(x)\om(y,z)- \phi_2(y)\om(x,z) & =
\om(y,x\cdot z)-\om (x,y\cdot z)\label{11}\\
a\cdot \om(y,z)+\phi_1 (y\cdot z)a & = \phi_2(y)\phi_1(z)a \label{12} \\
[\phi_2(x),\phi_2(y)] & = 0 \label{13} \\
\phi_2(y)(a\cdot c) & = a\cdot (\phi_2(y)c) \label{14} \\
a\cdot (\phi_1(z)b) & = b\cdot (\phi_1(z)a) \label{15}\\
\phi_1(z)\om(x,y)- \phi_1(y)\om(x,z) & =
\om(x\cdot z,y)-\om (x\cdot y,z)\label{16}\\
\om(x,y)\cdot c+\phi_2 (x\cdot y)c & = \phi_1(y)\phi_2(x)c \label{17} \\
[\phi_1(x),\phi_1(y)] & = 0 \label{18} \\
\phi_1(z)(a\cdot b) & = (\phi_1(z)a)\cdot b \label{19} \\
(\phi_2(x)c)\cdot b & = (\phi_2(x)b)\cdot c \label{20}
\end{align}
for all $a,b,c \in \La$ and  $x,y,z \in \Lb$.
\end{prop}

\begin{proof}
Let $u=(a,x), v=(b,y), w=(c,z)$ denote three arbitrary elements of $\Lg$.
Let us first consider the equation
\eqref{lr5} for the product, i.e., $[u,v] =u\kringel v - v\kringel
u$. Using \eqref{lie-algebra}, \eqref{produkt} and the commutativity
of the LR-product in $\La$ we obtain
\begin{align*}
[u,v] & =(\phi(x)b-\phi(y)a+\Om(x,y),[x,y])\\
u\kringel v -v\kringel u& =((\phi_2(x)-\phi_1(x))b-(\phi_2(y)-\phi_1(y))a
+\om(x,y)-\om(y,x),[x,y]).
\end{align*}
Suppose that the two expressions are equal for all $a,b\in\La$ and
$x,y\in \Lb$. For $a=b=0$ we obtain $\om(x,y)-\om(y,x)=\Om(x,y)$. Taking this into
account, $a=0$ implies $\phi_2(x)-\phi_1(x)= \phi(x)$. 
Conversely, equations \eqref{9} and \eqref{10} imply \eqref{lr5}.\\
A similar computations shows that \eqref{lr1} corresponds to the equations
$\ref{11},\ldots , \ref{15}$, and \eqref{lr2} corresponds to
$\ref{16},\ldots , \ref{20}$.
\end{proof}

\begin{cor}
Assume that $\Lg= \La \rtimes_{\phi} \Lb$ is a semidirect product of an abelian
Lie algebra $\La$ and a Lie algebra $\Lb$ by a representation
$\phi\colon \Lb \ra \End(\La)=\Der(\La)$, i.e., we have a split exact sequence
\begin{equation*}
0 \rightarrow \La \xrightarrow{\iota} \Lg \xrightarrow{\pi}
\Lb \rightarrow 0.
\end{equation*}
If $\Lb$ admits an LR-structure $(x,y)\mapsto x\cdot y$ such that
$\phi(x\cdot y)=0$ for all $x,y\in \Lb$, then also $\Lg$ admits an LR-structure.
\end{cor}

\begin{proof}
Because the short exact sequence is split, the $2$-cocycle $\Om$ in the Lie bracket
of $\Lg$ is trivial, i.e., $\Om(x,y)=0$. Let $a\cdot b=0$ be the trivial product
on $\La$ and take $\phi_1=0$, $\phi_2=\phi$ and $\om(x,y)=0$. Assume that $(x,y)\mapsto x\cdot y$ is
an LR-product. Then all conditions of Proposition $\ref{phi12}$ are satisfied except for
\eqref{17} and \eqref{13} requiring
\begin{align*}
\phi(x\cdot y) & = 0\\
[\phi(x),\phi(y)] & =0.
\end{align*}
But we have \eqref{17} by assumption, and since $\phi$ is a representation it follows
\[
0=\phi(x\cdot y-y\cdot x)=\phi([x,y])=[\phi(x),\phi(y)].
\]
Hence \eqref{produkt} defines an LR-structure on $\Lg$, given by
\[
(a,x)\kringel (b,y)=(\phi(x)b, x\cdot y).
\]
\end{proof}

\begin{cor}\label{trivial}
Suppose that the LR-products on $\La$ and $\Lb$ are trivial.
Hence $\Lb$ is also abelian.
Then \eqref{produkt} defines an LR-structure
on $\Lg$ if and only if the following conditions hold:

\begin{align}
\om(x,y)-\om(y,x) & = \Om(x,y)\\
\phi_2(x)-\phi_1(x) & = \phi(x)\\
\phi_2(x)\om(y,z) & = \phi_2(y)\om(x,z)\\
\phi_2(x)\phi_1(y) & = 0\\
[\phi_2(x),\phi_2(y)] & = 0\\
\phi_1(z)\om(x,y) & = \phi_1(y)\om(x,z)\\
\phi_1(x)\phi_2(y) & = 0\\
[\phi_1(x),\phi_1(y)] & = 0
\end{align}
\end{cor}

We can apply this corollary as follows:

\begin{prop}\label{iso}
Let $\Lg=(\La,\Lb,\phi,\Om)$ be an extension in which both $\La$ and $\Lb$ are abelian.
If there exists an $e\in \Lb$ such that $\phi(e)\in\End(\La)$ is an
isomorphism, then $\Lg$ admits an LR-structure. In fact, in that
case \eqref{produkt} defines an LR-product, where $\phi_1=0,\, \phi_2=\phi$,
the product on $\La$ and $\Lb$ is trivial, and
\begin{align*}
\om(x,y) & = \phi(e)^{-1}\phi(x)\Om(e,y).
\end{align*}
\end{prop}

\begin{proof}
We have to show that the above conditions of corollary
$\ref{trivial}$ are satisfied. Applying $\phi(e)^{-1}$ to \eqref{6}
with $z=e$ it follows
$\Om(x,y)-\phi(e)^{-1}\phi(x)\Om(e,y)+\phi(e)^{-1}\phi(y)\Om(e,x)=0$.
This just means that $\Om(x,y)=\om(x,y)-\om(y,x)$. Furthermore we
have, since $\phi(x)\phi(y)=\phi(y)\phi(x)$ for all $x,y \in \Lb$,
\begin{align*}
\phi(x)\om(y,z)-\phi(y)\om(x,z) & = \phi(x)\phi(e)^{-1}\phi(y)\Om(e,z)-
\phi(y)\phi(e)^{-1}\phi(x)\Om(e,z) \\
 & = \phi(e)^{-1}(\phi(x)\phi(y)-\phi(y)\phi(x))\Om(e,z)\\
 & = 0.
\end{align*}
All the other conditions follow trivially from $\phi_1=0$. Hence the product defines
an LR-structure.
\end{proof}

\end{document}